\documentclass [12pt,a4paper]{amsart} 
\usepackage{graphicx}
\usepackage[utf8]{inputenc} 
\usepackage[top=2.5cm,bottom=3.0cm,left=2.5cm,right=2.5cm]{geometry}
\usepackage{amssymb,amsmath,amsfonts,amsthm,mathtools,mathrsfs}
\usepackage{xcolor}
\usepackage{helvet}
\usepackage{hyperref}
\usepackage{cleveref}
\newtheorem{theorem}{Theorem}[section]
\newtheorem{lemma}[theorem]{Lemma}
\newtheorem{proposition}[theorem]{Proposition}
\newtheorem{corollary}[theorem]{Corollary}

\newtheorem{theoremletra}{Theorem}

\newtheorem*{prop1}{Proposition A.1}
\newtheorem*{prop2}{Proposition A.2}
\theoremstyle{definition}
\newtheorem{definition}{Definition}

\newtheorem*{question}{Question}
\theoremstyle{remark}

\newtheorem*{notation*}{Notation}
\usepackage[T1]{fontenc}
\usepackage[color=red!40]{todonotes}
\newcommand{\HH}{\mathbb{H}}

\newcommand{\GL}{\mathrm{GL}}
\newcommand{\SL}{\mathrm{SL}}
\newcommand{\PSL}{\mathrm{PSL}}
\newcommand{\OO}{\mathcal{O}}
\newcommand{\Spin}{\mathrm{Spin}}
\newcommand{\SO}{\mathrm{SO}}
\newcommand{\Br}{\mathrm{Br}}
\newcommand{\bs}{\backslash}

\newcommand{\Z}{\mathbb{Z}}
\newcommand{\Q}{\mathbb{Q}}
\newcommand{\R}{\mathbb{R}}
\newcommand{\C}{\mathbb{C}}
\newcommand{\gd}{\mathcal{O}_{d}}
\newcommand{\N}{\mathrm{N}}

\newcommand{\tr}{\mathrm{tr}}

\newcommand{\Od}{\mathcal{O}_d}

\newcommand{\Isom}{\mathrm{Isom}}
\newcommand{\Sal}{\mathrm{Sal}^{sq}}

\newcommand{\rank}{\mathrm{rank}}
\newcommand{\inv}{^{-1}}
\newcommand{\e}{e^{i\theta}}
\newcommand{\f}{e^{-i\theta}}
\newcommand{\SQ}{\Q^{\times}/\Q^{\times2}}

\newcommand{\Ll}{L_{\lambda}}
\usepackage{enumitem}
\newcommand{\Ram}{\mathrm{Ram}}

\begin{document}

\title[On Salem numbers of degree 4 and arithmetic hyperbolic orbifolds]
        {On Salem numbers of degree 4 and arithmetic hyperbolic orbifolds}

\author{Cayo D\'oria}
\author{Plinio G. P. Murillo}
 
\begin{abstract}

In this article, we construct an arithmetic hyperbolic $6-$orbifold $\mathcal{O}$ such that, any square-rootable Salem number of degree at most $4$ over $\Q$ is realized as the exponential of the length of a closed geodesic in $\mathcal{O}$. We also prove that $n=6$ is the minimal dimension among arithmetic hyperbolic orbifolds of the first type where it can be obtained. In an appendix, we establish a general relation between the discriminant of a Salem number and the determinant of a quadratic space which realizes it. In particular, for any $m,d>0$ we present a geometric proof of the existence of Salem numbers of degree $2m$ with discriminant $(-1)^{m+1}d$ in $\SQ$.

\end{abstract}

\maketitle 

\section{Introduction}
A real algebraic integer $\lambda > 1$ is a \emph{Salem number}\footnote{Since positive units of totally real number fields are Pisot numbers, Salem's original definition requires that a T number (which nowadays is called Salem number) has at least one non-real conjugate.} if its Galois conjugates outside the unit circle are precisely \(\lambda\) and $\lambda^{-1}$. In particular, \(\lambda\) has even degree over $\Q$. Although introduced by Salem motivated by problems in Fourier analysis (see \cite[Chapter III, Section 3]{Salem63}), these numbers have become relevant in several other branches of mathematics, such as algebraic geometry, dynamical systems and hyperbolic geometry (see \cite{GH01} and \cite{Smyth15}).

For instance, it is known for the last three decades that Salem numbers are related to the exponential length of closed geodesics of arithmetic hyperbolic \(2-\) and \(3-\)dimensional orbifolds (\cite{NR92},\cite{Sury92}). 

In 2019, Emery, Ratcliffe and Tschantz extended this connection to higher dimensions. 
\begin{theorem}[\cite{ERT19}, Thm 5.2 and Thm. 6.3]\label{th:ERT 5.2 e 6.3}
 Let $\lambda$ be an algebraic integer, $K$ be a totally real number field, and $n\in\Z$, with \(n \ge 2\). The following sentences are equivalent.
\begin{enumerate}
    \item $\lambda$ is a Salem number, $K\subset\Q(\lambda+\lambda^{-1})$ with $[\Q(\lambda):K]\leq n+1$;
    \item There exists a classical arithmetic lattice $\Gamma\subset\Isom(\HH^{n})$  of the first type, defined over $K$, and $\gamma\in\Gamma$ loxodromic with translation length \(\ell(\gamma)\), such that $\lambda=e^{\ell(\gamma)}$.
\end{enumerate}
\end{theorem}

We refer the reader to \Cref{sec:prel} for the relevant definitions concerning the \Cref{th:ERT 5.2 e 6.3}. Following \cite{CP25}, we say that a Salem number $\lambda$ is \textit{realized} by a subgroup $\Gamma<\Isom(\HH^n)$ if there exists $\gamma\in\Gamma$ loxodromic such that $\lambda=e^{\ell(\gamma)}$. In the case that \(\Gamma\) is discrete, we also say that the orbifold $\OO=\Gamma\bs\HH^n$ realizes $\lambda$. 

\Cref{th:ERT 5.2 e 6.3} implies that a classical arithmetic lattice of the first type in $\Isom(\HH^n)$ defined over $\Q$ realizes a Salem number $\lambda$ only if $[\Q(\lambda):\Q]\leq n+1$. Therefore, it is possible that a single orbifold realizes \textit{all} Salem numbers of a given degree. This is the main motivation of this article, and we summarize it in the following question.

\begin{question}
\hypertarget{q:question}{}
\textit{Is there a classical arithmetic hyperbolic lattice of the first type $\Gamma<\Isom(\HH^n)$ defined over $\Q$, with $n\geq 2m-1$, realizing all Salem numbers of degree less than or equal to $2m$? }
\end{question}

This question is not trivial even for \(m=1\). We prove that there does not exist a Fuchsian group realizing all Salem numbers of degree 2 (see \Cref{sec:2Salem}). However, the modular group $\Gamma=\SL_2(\Z)$ realizes $\lambda^2$ for any Salem number \(\lambda\) of degree 2 (see \cite[Corollary 1.5]{Sarnak82}). In this work we obtain a generalization of this result. 

Indeed, we are able to prove the existence of a nonuniform \(6\)-dimensional arithmetic lattice that realizes the squares of all Salem numbers of degree at most \(4\). Furthermore, we prove that this orbifold realizes a larger class of Salem numbers, which we discuss below.

One of the main contributions of \cite{ERT19} is the introduction of the so-called \textit{square-rootable Salem numbers} (see \Cref{def:square-rootable}). These particular numbers appear as squares of exponential length of closed geodesics in arithmetic hyperbolic orbifolds of the first type in odd dimensions. It follows from \cite[Lem. 7.4]{ERT19} that the square of any Salem number is square-rootable of the same degree, however, not any square-rootable Salem number is given in this way (see \Cref{sec:2Salem}). Some properties and applications of these special numbers were studied in \cite{BLMT22}, \cite{MMRT25} and \cite{G24}. 

Now we can state our first result, where we give a positive answer to the \hyperlink{q:question}{Question} for Salem numbers of degree less than or equal to $4$ which are square-rootable over $\Q$. In what follows, for simplicity, we will use the notation
\begin{align*}
    \Sal_m = \{\text{Salem numbers of degree } 2m, \text{ square-rootable over } \Q\}. 
\end{align*}

\begin{theoremletra}\label{th:dim 6}
    There is a non-cocompact arithmetic hyperbolic \(6\)-orbifold \(\mathcal{O}=\Gamma \backslash \HH^6\) that realizes all Salem number of degree 2, and all elements of $\Sal_2.$ In particular, for any Salem number $\lambda$ of degree $4$, there exists a closed geodesic in \(\mathcal{O}\) with length \(2\log(\lambda)\).
\end{theoremletra}

The orbifold $\mathcal{O}$ is constructed in \Cref{sec:dim 6}, and it is a generalization of the modular surface if we consider the Vahlen group as a model for the group of isometries of $\HH^6$, which is formed by \(2 \times 2\) matrices with entries in a certain Clifford algebra. We give an overview about this model in \Cref{sec:vahlen group}. The arithmetic part of the construction uses special features of Salem numbers of degree \(4\) which are presented in \Cref{sec:4Salem}, and is also an application of the celebrated Lagrange Four Square Theorem.

We also prove that $n=6$ is the minimal dimension of an orbifold where all elements of \(\Sal_2\) can be realized by a fixed orbifold. In fact, we prove a stronger result. We say that a commensurability class realizes a Salem number if some lattice in the class realizes it.

\begin{theoremletra}\label{th:3-5}
    There is no a commensurability class of arithmetic lattices of first type in \(\Isom(\HH^n)\), $3\leq n<6$, realizing all Salem numbers in \(\Sal_2\).
\end{theoremletra}

Commensurability classes of arithmetic lattices of the first type in \(\Isom(\HH^n)\) are determined by orthogonal groups $\SO_q(K)$, where $q$ is an admissible quadratic form over a totally real number field $K$ (see \Cref{sec:aho} for the admissibility condition). For the purpose of this article it is enough to consider $K=\Q$ (\Cref{prop:rational_base_field}). The proof of \Cref{th:3-5} uses algebraic invariants of quadratic forms to constrain the Salem numbers that are realized by a given commensurability class. The case \(n=3\) uses the determinant. For $n=4,5$ we use the Witt and Hasse invariants, combined with local properties of quaternions algebras and number fields. Since the invariants become more complex as the dimension grows, we split the proof into three theorems, one for each dimension, which are the content of \Cref{sec:dim 4-5}.

Let us end this introduction by commenting some results that we have added in an Appendix. The key part of the proof of \Cref{th:3-5} for $n=3$ relies on a general connection between the characteristic polynomial of elements in $\SO_q(\Q)$, and the determinant of $q$. Let $F_\lambda\in\Z[x]$ be the minimal polynomial of a Salem number $\lambda$, and suppose that there is an isometry of $q$ with characteristic polynomial $F_\lambda$. By \cite{BF15}, this implies that $F_\lambda(1)F_\lambda(-1)=\det(q)$ in $\SQ$. This imposes a restriction for a Salem number to be realized by an arithmetic subgroup of $\Isom(\HH^n)$ of the first type. The class of $F_\lambda(1)F_\lambda(-1)$ in \(\SQ\) has also appeared in recent works, for instance \cite{BF23}, \cite{BT20}, and \cite{GM02}. It turns out that the class of $F_\lambda(1)F_\lambda(-1)$ in $\SQ$ is related to the class of the discriminant $d_\lambda$ of $\lambda$. Although this connection is well-known to experts (e.g. \cite[p.85]{D16}), it is not formulated as a separated proposition. We have decided to provide an explicit statement and a proof of the following.
\begin{prop1}\label{prop:equality_discriminants}
Let $m\geq 2$ and $\lambda$ be a Salem number of degree \(2m\)  with minimal polynomial \(F_\lambda\) and discriminant $d_\lambda$. Then
    $d_\lambda=(-1)^{m}F_\lambda(1)F_{\lambda}(-1)$ in $\Q^{\times}/\Q^{\times 2}$.
\end{prop1}

As a consequence, we give a geometrical proof to the following number-theoretical result.
\begin{prop2}\label{prop:existence of salem with discriminant given}
  Given a square-free integer $d>0$ and an integer \(m>0\), there exists a Salem number  $\lambda$ of degree \(2m\)  with $d_\lambda \equiv (-1)^{m+1}d \mod \SQ$.   Furthermore,  \(\lambda\) can be chosen square-rootable over \(\Q\).
\end{prop2}

\section{Preliminaries}\label{sec:prel}

\subsection{Hyperbolic Manifolds}

The hyperboloid model of the hyperbolic $n$-space is given by
$$\mathbb{H}^n=\{x \in \mathbb{R}^{n+1};\  x_0^2-x_1^2-\cdots - x_n^2=1, \ x_0>0\},$$
with the metric $ds^2=-dx_0^2+dx_1^2+\cdots+dx_{n}^2$. In this way, $\HH^n$ is the unique complete, simply connected $n$-dimensional Riemannian manifold with constant sectional curvature equal to $-1$, up to isometry. The group of isometries $\Isom(\HH^n)$ of $\HH^n$ is isomorphic to the Lie group
\begin{align*}
\mathrm{O}(1,n)^+=\{ A \in \GL_{n+1}(\R) \:|\: A^tJA=J, \mbox{ and } A \mbox{ preserves } \mathbb{H}^n\}  
\end{align*}
where $J$ is the rank $n+1$ diagonal matrix $J=diag(1,-1,-1,\ldots,-1)$. The identity component $\mathrm{SO}^{\circ}(1,n)$ of $\SO(1,n)=\mathrm{O}^+(1,n) \cap\SL_{n+1}(\R)$ is isomorphic to the group $\mathrm{Isom}^{+}(\mathbb{H}^n)$ of orientation preserving isometries of $\mathbb{H}^{n}$. Given a lattice $\Gamma\subset\mathrm{Isom}(\mathbb{H}^n)$, i.e, a discrete subgroup having finite covolume with respect to the Haar measure of $\Isom(\HH^n)$, the associated quotient space $M=\Gamma\backslash\mathbb{H}^n$ is a finite volume  {\it hyperbolic orbifold}, which is a manifold whenever $\Gamma$ is torsion-free. 

We recall that an element $\gamma\in\mathrm{Isom}(\mathbb{H}^n)$ is called
\begin{itemize}
\item {\it elliptic}, if it has a fixed point on $\mathbb{H}^n$.
\item {\it parabolic}, if it has exactly one fixed point in $\partial \mathbb{H}^n$.
\item {\it loxodromic}, if it has exactly  two fixed points in $\partial \mathbb{H}^n$.
\end{itemize}

The {\it translation length} of a loxodromic element $\gamma\in\Isom(\HH^{n})$ is defined by
$$\ell(\gamma):= \inf\limits_{x\in \mathbb{H}^n} d_{\HH^n}(\gamma \cdot x,x).$$

We recall that the eigenvalues of $\gamma$ out of the unit circle are precisely $e^{\ell(\gamma)}$ and $e^{-\ell(\gamma)}$ (see \cite[Proposition 1]{Greenberg62}). Therefore, a subgroup $\Gamma\subset\Isom(\HH^n)$ realizes a Salem number $\lambda$ if and only if there is $\gamma\in\Gamma$ loxodromic with an eigenvalue equal to $\lambda$.

\subsection{Arithmetic hyperbolic orbifolds and square-rootable Salem numbers}\label{sec:aho}
Let $K$ be a totally real number field, and let $q$ be a quadratic form defined over $K$. We say that $q$ is \textit{admissible} if $q$ has signature $(1,n)$ over $\mathbb{R}$, and for any non-trivial embedding $\sigma: K\rightarrow\mathbb{R}$ the quadratic form $f^{\sigma}$ is definite.

The conditions for $q$ being admissible imply that $\SO_{q}$ is a semisimple algebraic $K$-group such that $\SO_{q}(\mathbb{R})\simeq \SO(1,n)$ and $\SO_{q^{\sigma}}(\R)$ is compact for any non-trivial embedding $\sigma: K \rightarrow\mathbb{R}$. By Weil's restriction of scalars, $\SO_{q}(\mathcal{O}_{K})$ embeds as an arithmetic subgroup of $\SO(1,n)$. 

We recall that two subgroups $\Gamma_1,\Gamma_2$ of $\Isom(\HH^n)$ are \textit{commensurable} if there is $g\in\Isom(\HH^n)$ such that $g\Gamma_1g\inv\cap\Gamma_2$ has finite index in both $\Gamma_1$ and $\Gamma_2$. A subgroup $\Gamma\subset\Isom(\mathbb{H}^{n})$ commensurable with $\SO_{q}(\mathcal{O}_{K})$ for some admissible $q$ is called an \textit{arithmetic subgroup of $\Isom(\mathbb{H}^{n})$ of the first type defined over} $K$. The field $K$ is a commensurability invariant of $\Gamma$ (see \cite[Section 3]{M11}). The quotient $M=\Gamma\backslash\mathbb{H}^{n}$ is called an \textit{arithmetic hyperbolic orbifold of the first type defined over} $K$. If $\Gamma$ is torsion-free,  $M=\Gamma\backslash\mathbb{H}^{n}$ is an \textit{arithmetic hyperbolic manifold of the first type defined over} $K$.
Moreover, $\Gamma$ is said to be \textit{classical} if $\Gamma<\SO_q(K)$ for some admissible quadratic form $q$. For $n$ even, any arithmetic subgroup of $\Isom(\HH^{n})$ is classical \cite[Lemma 4.2]{ERT19}. For $n$ odd, this is not longer true, but the group $\Gamma^{(2)}$ generated by the squares of elements of $\Gamma$ is always classical \cite[Lemma 4.5]{ERT19}.

As we mentioned in the Introduction, Emery, Ratcliffe and Tshantz proved that Salem numbers are precisely the exponential of the translation length of loxodromic elements in classical arithmetic subgroups of $\Isom(\HH^n)$ of the first type (see \Cref{th:ERT 5.2 e 6.3}). They also introduced a special class of Salem numbers, whose definition we recall below.

\begin{definition}\label{def:square-rootable}
Let $\lambda$ be a Salem number, $K$ be a number field contained in $\mathbb{Q}(\lambda+\lambda^{-1})$ and let $f(x)$ be the minimal polynomial of $\lambda$ over $K$. We say that $\lambda$ is \textit{square-rootable over $K$} if there is a totally positive $\alpha\in K$ and a monic palindromic polynomial $g(x)$, whose even degree coefficients are in $K$ and whose odd degree coefficients are in $\sqrt{\alpha}K$, such that $f(x^2)=g(x)g(-x)$.
\end{definition}

\begin{theorem}[\cite{ERT19}, Thm. 7.6 and Thm. 7.7]
 Let $\lambda$ be a Salem number, $K$ be a totally real number field, and $n\in\Z$ odd. Then the following sentences are equivalent.

 \begin{enumerate}
    \item $K$ is a subfield of $\Q(\lambda+\lambda^{-1})$, $\lambda$ is square-rootable Salem number over $K$, and $[\Q(\lambda):K]\leq n+1$.
    \item There exists $\Gamma\subset\Isom(\HH^{n})$ an arithmetic lattice of the first type defined over $K$ and $\gamma\in\Gamma$ loxodromic such that $\lambda^{1/2}=e^{\ell(\gamma)}$.
\end{enumerate}
\end{theorem}

\subsection{Clifford algebras and the Spin group}\label{sec:Clifford algebras}

In this section we will recall the construction of Clifford algebras and spin groups. For further details, we refer the reader to \cite[Section 2]{Mur19} and the references therein.
 
 Let $K$ be a field with char~$K\neq 2$, and $q$ be a non-degenerate quadratic form on $K^n$ with associated bilinear form $\Phi$. We will denote the rank \(n\) of \(q\) by \(\rank(q)\). For an orthogonal basis $i_1,\dots,i_{n}$ of $K^n$ with respect to $q,$ the \textit{Clifford algebra of $q$} is the unitary associative algebra over $K$ generated by  $i_1,\dots,i_{n}$ with relations $i_{\nu}^2=q(i_{\nu})$ and $i_{\nu}i_{\mu}=-i_{\mu}i_{\nu}$ for $\mu,\nu\in\{1,\ldots,n\}, \ \mu \neq \nu.$ Let $\mathcal{P}_{n}$ be the power set of $\{1,...,n\}$. For $M=\{\mu_1,\dots,\mu_\kappa\}\in \mathcal{P}_{n}$ with $\mu_1<\dots<\mu_\kappa,$ we define $i_M=i_{\mu_1}\cdot...\cdot i_{\mu_\kappa}$, where we adopt the convention $e_{\emptyset}=1$. The $2^n$ elements $i_M$, $M\in\mathcal{P}_n$ determine a basis for $\mathcal{C}(q)$, so every element in $\mathcal{C}(q)$ can be written uniquely as  $s=\sum_{M \in \mathcal{P}_{n}}s_{M}i_{M}$, with $s_M\in K$. We identify $K$ with $K\cdot 1$, and $K^n$ with the $K$-linear subspace in $\mathcal{C}(q)$ generated by $i_1, i_2,\ldots, i_n$. The \textit{vectors} on $\mathcal{C}(q)$ are the  elements in $K\cdot 1\oplus K^n$, that is, the $K$-linear combinations of $1,i_1,\ldots,i_n$.
 
The algebra $\mathcal{C}(q)$ has an anti-involution $*$, which acts on the elements $i_M$ by $i_{M}^{*}=(-1)^{\kappa(\kappa-1)/2}\ i_{M}$, where \(\kappa=|M|\). 
The span of the elements $i_M$ with $|M|$ even is a subalgebra $\mathcal{C}^{+}(q)$ of $\mathcal{C}(q),$ called the \textit{even Clifford subalgebra of $q$}. 

The spin group of $q$ is defined as
\begin{equation*}
\Spin_{q}(K):=\left\{ s \in  \mathcal{C}^{+}(q)\ \Big| \ sK^ns^*\subseteq K^n, \ ss^*=1 \right\}.    
\end{equation*}

For an element $s\in \Spin_{q}$, the map $\varphi_{s}:K^n \to K^n$ given by $\varphi_{s}(x)=sxs^{-1}$ preserves the quadratic form $q$ and defines a homomorphism of algebraic $K$-groups
\begin{align} 
\varphi: \Spin_{q}&\rightarrow \SO_{q}, \label{phi}  \\
s & \mapsto\varphi_{s}. \nonumber
\end{align}

The kernel of $\varphi$ is the set $\lbrace1,-1\rbrace$ and, if $q$ is isotropic (i.e. $q(v)=0$ for some non-zero $v\in K^n$), then
\begin{equation}\label{eq:index spin vs so}
\SO_{q}(K)/\varphi(\Spin_{q}(K))\simeq K^{\times}/(K^{\times})^{2},    
\end{equation}
see \cite[II.2.3, II.2.6, II.3.3 and II.3.7]{Chevalley54} for more details. For future reference we will denote by \(q_n\) the following rank \(n+1\) quadratic form
\begin{equation}\label{eq: canonical form}
    q_n=x_0^2-x_1^2-\cdots-x_n^2
\end{equation}
over any field. In the case $K=\mathbb{R}$, the corresponding spin group is denoted by $\Spin(1,n)$. 

\subsection{The Vahlen group}\label{sec:vahlen group} The group of orientation-preserving isometries of $\HH^n$ can be represented by $2\times2$ matrices with entries in a Clifford algebra, the so-called \textit{Vahlen's group}. This generalizes the groups \(\PSL_2(\R),\PSL_2(\C)\) as isometry groups of hyperbolic spaces for $n=2$ and $n=3$. In this section we will recall the definition of the Vahlen group and its relation with $\HH^n$. For further details, we refer the reader to \cite{Ahlfors85}, \cite{EGM87}, \cite{MWW89} and \cite{Wat93}, and the references therein.

Let $K$ be a field. For each \(m \in \Z_{>0}\) consider the quadratic form of rank $m$ over $K$ given by $h_m=-x_1^2-x_2^2-\cdots-x_m^2$, and denote by  \(\mathcal{C}_{m}(K)\) the Clifford algebra  $\mathcal{C}(h_m)$. In order to keep the usual notation, we will denote by \(\{i_1,\ldots,i_m\}\) the generators of \(\mathcal{C}_{m}(K)\) .

For any vector \(\alpha=a_0+a_1i_1+\ldots+a_m i_m,\) we define the conjugated of $\alpha$ as \(\overline{\alpha}=a_0-a_1i_1-\ldots-a_mi_m\). The norm and trace of \(\alpha\) are defined respectively as
\begin{align*}
    \rm{Norm}(\alpha):= & ~ \alpha \overline{\alpha} = a_0^2+\ldots+a_m^2. \\
    \rm{tr}(\alpha):= & ~ \alpha + \overline{\alpha} =2a_0.
\end{align*}
In this way, any vector \(\alpha\in\mathcal{C}_{m}(K)\) satisfies the quadratic equation over \(K\)
\begin{align}\label{eq: quadratic equation for vector}
    \alpha^2-\tr(\alpha) \alpha+\rm{Norm}(\alpha)=0.
\end{align}

The\textit{ Clifford group} \(\Gamma_m(K)\) is the multiplicative subgroup of \(\mathcal{C}_{m}(K)\) generated by the nonzero vectors. For \(K=\R\) and \(m=0,1,2\) we have respectively \(\Gamma_{m} \cup \{0\}=\R,\C,\mathcal{H}\), where \(\mathcal{H}\) denotes the Hamilton's quaternion algebra.  When \(K=\R\) we will drop the field of definition and denote \(\mathcal{C}_{m}(\R)\) (respectively \(\Gamma_m(\R))\) only by \(\mathcal{C}_{m}\) (respectively \(\Gamma_m\)).

The \textit{Vahlen group}\footnote{See the paragraph after Definition 1 in \cite{Ahlfors85} for an explanation for this shortened definition.} is defined as follows
    \begin{align*}
        \SL_2(\mathcal{C}_{m})=\left\{ \begin{pmatrix}
            a & b \\ c & d
        \end{pmatrix} : a,b,c,d \in \Gamma_m \cup\{0\},~ ad^*-bc^*=1 \text{ and } ab^*, cd^* \in \R^{m+1}\right\}.
    \end{align*}

The conditions on the coefficients of the elements of \(\SL_2(\mathcal{C}_{m})\) guarantees that \(\SL_2(\mathcal{C}_{m})\) is a group, with respect to matrix multiplication. On the other hand, the set of vectors of the Clifford algebra $\mathcal{C}_{m+1}$ contains a copy of the \((m+2)\)-dimensional hyperbolic space given by

\[\HH^{m+2}=\{x_0+x_1i_1+\ldots+x_{m+1}i_{m+1} \in \mathcal{C}_{m+1} : x_{m+1}>0\},\]
Moreover, \(\SL_2(\mathcal{C}_{m})\) acts on  \(\HH^{m+2}\) as 
  \begin{align*}
      \begin{pmatrix}
            a & b \\ c & d
        \end{pmatrix} \cdot v=(av+b)(cv+d)\inv.
  \end{align*}
\begin{theorem}\cite[Theorem 5]{Wat93}
  The group \(\SL_2(\mathcal{C}_{m})\) acts isometrically on \(\HH^{m+2}\).  \\
  Moreover, the full group of orientation-preserving isometries of \(\HH^{m+2}\) is isomorphic to  \(\rm{PSL}_2(\mathcal{C}_{m})=\SL_2(\mathcal{C}_{m})/\{\pm \rm{I}\}.\)
\end{theorem}

For any subring \(R \subset \R\) we denote by \(\mathcal{C}_{m}(R)\) the  subring of \(\mathcal{C}_{m}\) whose elements have coefficients in \(R\). The subgroup \(\SL_2(\mathcal{C}_{m}(\Z))\) is a natural generalization of the well-known groups \(\SL_2(\Z)\) and the \(\SL_2(\Z[i])\). Indeed, we have the following.

\begin{theorem}\cite[Theorem 3]{MWW89}\label{th:arithmetic subgroup Vahlen}
For any $m\geq 0,$ \(\SL_2(\mathcal{C}_{m}(\Z))\) is an arithmetic subgroup of \(\SL_2(\mathcal{C}_{m})\).
\end{theorem}

We finish this section with an explicit isomorphism between \(\SL_2(\mathcal{C}_{n-2}(K))\) and \\ \(\Spin_{q_n}(K)\), when $K\subset\R$ is a number field, and \(q_n\) is as in \Cref{eq: canonical form}. 

Consider a canonical basis for the Clifford algebra \(\mathcal{C}(q_{n})\) written as \(f_0,f_1,f_2,i_1,\ldots,i_{n-2}\), that is, \(f_0^2=1, f_1^2=f_2^2=-1\) and \(i_j^2=-1\) for \(j=1,\ldots,n-2\). The notation of the basis has been chosen so that the elements \(i_1,i_1,\ldots,i_{n-2}\) generates the algebra \(\mathcal{C}_{n-2}(K)\) as a subalgebra of \(\mathcal{C}(q_n)\). Now, consider the following elements in \(\mathcal{C}(q_n)\):
\begin{align*}
    \tau_0=\frac{1}{2}(f_0+f_1), ~ & ~ \tau_1=\frac{1}{2}(f_0-f_1), \\
    u=\tau_1\tau_0, ~ & ~ w_1=\tau_1f_2, \\
    w_0 = \tau_0 f_2, ~ & ~ v=\tau_0\tau_1.
\end{align*}
For any subfield \(K \subset \R\), there exists an injective \(K\)-algebra homomorphism \(\iota\) from \(\mathcal{C}_{n-2}(K)\) to \(\mathcal{C}_{q_n}^+(K)\), whose image on the generator are given by 
\(\iota(i_j)=f_0f_1f_2i_j,
\text{ for } j=1,\ldots,n-2\)
(see \cite[Proposition 2.4]{EGM87}). With this immersion we can present the following explicit map between \(\SL_2(\mathcal{C}_{n-2}(K))\) and \(\Spin_{q_n}(K).\) 

\begin{theorem}\cite[Thm. 4.1]{EGM87} 
  For any subfield \(K \subset \R\), the map
    \begin{align}
    \psi: \SL_2(\mathcal{C}_{n-2}(K)) &\to \mathcal{C}^+_{q_n}(K) \label{map psi} \\ 
        \begin{pmatrix}
            a & b \\ c & d
        \end{pmatrix}& \mapsto \iota(a)u+\iota(b)w_1+\iota(c)w_0+\iota(d)v 
    \end{align}
 gives an \(K\)-isomorphism of \(K\)-algebraic groups between \(\SL_2(\mathcal{C}_{n-2}(K))\) and \(\Spin_{q_n}(K).\)   
\end{theorem}

\subsection{The Witt invariant}\label{sec:Witt} Let $q$
 be a quadratic form over $K$. The structure of $\mathcal{C}(q),$ depends on the parity of $\rank(q)$. Let $z=i_1i_2\cdots i_n$. If $\rank(q)$ is odd, $z$ lies in the center of $\mathcal{C}(q),$ and $\mathcal{C}^+(q)$ is a central simple algebra over $K$. When $\rank(q)$ is even, $z$ lies in the center of $\mathcal{C}^+(q),$ and $\mathcal{C}(q)$ is a central simple algebra over $K$. The \emph{Witt invariant} of $q$, denoted by $c(q)$, is the element in the Brauer group $\Br(K)$ of $K$ given by 
\[c(q)=\begin{cases}
    [\mathcal{C}^{+}(q)],& \text{  if $\rank(q)$ is odd}.\\
    [\mathcal{C}(q)], &\text{  if $\rank(q)$ is even}.\\
\end{cases}\]
A quaternion algebra $A$ over $K$ is a central simple $K$-algebra of dimension $4$. It has a standard basis, given by $1, i, j , ij$, where $i^2=c, j^2=d, ij=-ji, c, d\in K^{*}$, and $A$ can be represented by the Hilbert symbol $\left(\frac{c,d}K\right)$. The \textit{Hasse invariant}, $s(q)$, which is also an element in $\Br(K)$ of $K$ is given by the class of  \[\prod_{i<j}\left(\frac{a_i,a_j}K\right),\] where $\{a_1,a_2,\ldots,a_n\}$ is a diagonalization of $q$. The Hasse invariant is independent of the diagonalization, and is related to the Witt invariant in $\Br(K)$ as follows (see \cite[Prop. 3.20]{Lam05}).
\begin{proposition}\label{prop:hasse-witt} $s(q)$ and $c(q)$ satisfy the following conditions:
\begin{equation*}
    c(q)=
\begin{cases}
    s(q) &\text{ if } \rank(q)\equiv 1,2\mod 8.\\
    s(q)\cdot\left(\frac{-1, -d(q)}{K}\right) &\text{ if } \rank(q)\equiv 3,4\mod 8.\\ 
    s(q)\cdot\left(\frac{-1,-1}{K}\right)&\text{ if } \rank(q)\equiv 5,6\mod 8.\\
    s(q)\cdot\left(\frac{-1, d(q)}{K}\right)&\text{ if } \rank(q)\equiv 7,8\mod 8.\\
\end{cases}
\end{equation*}

where $d(q)$ denotes de determinant of $q$.
\end{proposition}

In general, as elements of $\Br(K)$, $c(q)$ and $s(q)$ can be expressed as products of quaternion algebras, and then both have order two \cite[Cor. V. 3.14]{Lam05}. In particular, when $K$ is a number field, both $c(q)$ and $s(q)$ can be represented by the class of some quaternion algebras \cite[Thm. VII. 2.6]{Milne11}.

Let $K$ be a number field, $A$ be a quaternion algebra over $K$, and $\nu$ be a place of $K$ with completion $K_\nu$. We recall that $A$ \textit{ramifies} at $\nu$ if $A_\nu=A\otimes_KK_\nu$ is a division algebra. Otherwise, $A$ \textit{splits} at $\nu$, in which case $A\cong M_2(K_\nu)$. By the classification of quaternion algebras over number fields, the set $\Ram(A)$ of places of $K$ where $A$ ramifies is always finite of even cardinality. The set $\Ram_f(A)$ denotes the finite places of $\Ram(A)$.

\begin{lemma}\label{lem:ram at inf}
    Let $(E,q)$ be a quadratic space over $\Q$,  and suppose that $q$ is admissible. Then
    \begin{align*}
    &\text{ If  }  \dim(E)\equiv 0,1,2,3\hspace{-3mm}\mod 8,\text{ then } c(q) \text{ splits at infinity.}\\
    &\text{ If } \dim(E)\equiv 4,5,6,7\hspace{-3mm}\mod 8,\text{ then } c(q) \text{ ramifies at infinity.} 
    \end{align*}
\end{lemma}

\begin{proof}
Let $n=\dim(E)$, $E_{\R}=E\otimes_{\Q}\R$ and denote by $q_\R$ to $q$ as a quadratic form over $\R$. By the admissibility condition we have that $(E_\R,q_\R)\cong(\R^{n},x_0^2-x_1^2-\cdots-x_n^2)$, therefore $d(q_\R)=(-1)^{n}$ and
\begin{equation*}
    s(q_\R)=
\begin{cases}
    M_2(\R) &\text{ if } \dim(E)\equiv 1,2, 5,6\mod 8.\\
    \left(\frac{-1,-1}{\R}\right) &\text{ if } \dim(E)\equiv 3,4,7,8\mod 8.\\
    \end{cases}
\end{equation*}
The result now follows from \Cref{prop:hasse-witt}.  
\end{proof}

\section{Salem numbers of degree at most 4}\label{sec:4Salem}

In this work, the special case of Salem numbers of degrees two and four play a central role. In what follows, we will explore some features of these numbers.

\subsection{Salem numbers of degree 2}\label{sec:2Salem}

Let \(\lambda\) be a Salem number of degree 2 with minimal polynomial \(F_{\lambda}(x)=x^2-Nx+1,\) where \(N\in\mathbb{N}\) and \(N^2-4>0\).

We note that 
\begin{align}
    F_\lambda(x^2)= & (x^2-\sqrt{N+2}~ x+1)(x^2+\sqrt{N+2}~ x+1). \label{eq:2-Salem are SR} 
\end{align}

Thus, \(\lambda\) is a square-rootable Salem number over \(\Q\). In particular, any totally positive fundamental unit of a real quadratic number field is square-rootable Salem number over $\Q$ which is not a square.

\begin{proposition}
There is no an arithmetic hyperbolic \(2\)-orbifold which realizes all Salem numbers of degree \(2.\)
\end{proposition}
\begin{proof}
Suppose that some arithmetic \(2-\)orbifold \(\mathcal{O}=\Gamma \backslash \HH^2\) realizes all Salem numbers of degree \(2\) for some lattice \(\Gamma < \PSL_2(\R)\). Then for any Salem number \(\lambda\) of degree \(2\) there exists \(\ell\) in the length spectrum of \(\mathcal{O}\) with \(\lambda=e^\ell\). Hence, we can find a hyperbolic element \(\gamma \in \Gamma\) with eigenvalues \(e^{\pm \frac{\ell}2}\), which implies in particular that \(\tr(\gamma)=e^{\frac{\ell}2}+e^{-\frac{\ell}2}\). Thus, the field \(k=\Q(\{\tr(\gamma) : \gamma \in \Gamma\})\) contains \(\lambda^{\frac{1}{2}} + \lambda^{-\frac 12}\) for all Salem number \(\lambda\) of degree \(2\). Moreover, \(k\) is a number field since all traces of \(\Gamma\) are algebraic and \(\Gamma\) is finitely generated (see \cite[Lemma 3.5.1]{MR02}).

On the other hand, if \(\lambda+\lambda\inv=N>2\), then \(\lambda^{\frac{1}{2}} + \lambda^{-\frac 12}=\sqrt{N+2}\). Hence, \(k\) contains infinitely many real quadratic extension of \(\Q,\) which gives a contradiction with the fact that \(k\) is a number field.    
\end{proof}

\subsection{Salem numbers of degree 4}
Let  \(\lambda\) be  a Salem number of degree \(4\) with conjugates \(\lambda\inv,\e,\f\) (with \(\theta \in (0,\pi))\) and consider \(F_\lambda=X^4+aX^3+bX^2+aX+1\) its minimal polynomial. We have an explicit simple expression for \(F_\lambda(1)F_\lambda(-1).\) Indeed, 
\begin{align*}
    F_\lambda(1)F_\lambda(-1) =~ & [(b+2)+2a]\cdot [(b+2)-2a] \\
    = ~  & (b+2)^2-4a^2.
\end{align*}

We associate to any Salem number \(\lambda\) of degree 4 the complex number \(\tau=\tau(\lambda)=\lambda\e+\lambda\inv\f\). We note that \(\tau\) is an imaginary quadratic algebraic integer since \(\tau\) is not real and satisfies
\begin{equation} \label{equation: tau}
 \tau^2 - (b-2)\tau + (a^2-2b)=0.  
\end{equation}

Indeed, since
\begin{equation}\label{eq:a and b}
\begin{split}
   a = & ~ -(\lambda+\lambda\inv+\e+\f).\\
   b = & ~ 2+(\lambda+\lambda\inv)(\e+\f). 
\end{split}
\end{equation}

and  \(\overline{\tau}=\lambda\inv\e+\lambda\f\), then \(\tau+\overline{\tau}=b-2\) and
\begin{equation}\label{norm of tau}
    \tau\overline{\tau} =  \lambda^2+\lambda^{-2}+e^{2i\theta}+e^{-2i\theta} = a^2-2b.
\end{equation}

It follows that the discriminant  \(\Delta=(b-2)^2-4(a^2-2b)\)  of \Cref{equation: tau} is negative.
We can rewrite \(\Delta\) as
\begin{equation}
\Delta = (b+2)^2-4a^2.
\end{equation}

Hence, the discriminant of the minimal polynomial of \(\tau\), over \(\Q\), is equal to \(F_\lambda(1)F_\lambda(-1)\). By \Cref{thm:discriminant vs F(1)F(-1)} we conclude that \(\Delta \equiv d_\lambda \mod \SQ\). For future reference, we summarize these observations into the following proposition.

\begin{proposition}\label{pro:field of tau}
For any Salem number \(\lambda\) of degree 4, \(\tau(\lambda)\) is an imaginary quadratic integer. Moreover, \(\Q(\tau(\lambda))=\Q(\sqrt{d_\lambda})\).   
\end{proposition}

\begin{lemma}\label{lem:min poly beta}
  Let \(\ell>0\). If $e^{\ell}$ is a Salem number of degree 4 with complex conjugates \(e^{\pm i\theta}\), then the minimal polynomial of $\beta=e^{\ell+i\theta}$ over $\Q$ is given by \[ r(x)=x^4-(4\cosh{\ell}\cos{\theta})x^3+(2+4\cosh^2{\ell}-4\sin^2{\theta})x^2-(4\cosh{\ell}\cos{\theta})x+1.\]
\end{lemma}

\begin{proof}
We have seen that $\tau=\tau(\lambda)$ is quadratic over $\Q$. On the other hand, $\beta$ is a root of $x^2-\tau x+1$ and \(\tau^2-4 \notin \Z,\) so $\beta$ has degree $4$ over $\Q$. Consider the polynomial 
\begin{align*}
    p(x)&=(x^2-\tau x+1)(x^2-\tau' x+1)\\
    &=x^4-(\tau+\tau')x^3 +(2-\tau\tau')x^2-(\tau+\tau')x+1\\
    &=x^4-(b-2)x^3+(2-a^2+2b)x^2-(b-2)x+1
\end{align*}

which lies in $\Z[x]$ and has $\beta$ as a root. Then, $p(x)$ is the minimal polynomial of $\beta$ over $\Q$. By \eqref{eq:a and b} we get that $a=-2(\cosh(\ell)+\cos(\theta))$ and $b=2+4\cosh\ell\cos\theta,$ so a direct computation shows that $p(x)=r(x).$
\end{proof}

\subsection{Square-rootable Salem numbers of degree 4}

When \(\lambda \in \Sal_2\), that is, \(\lambda\) is a  Salem number of degree \(4\) and square-rootable over \(\Q\), additionally to \(\tau(\lambda)\) we also attach to \(\lambda\) the complex number \[\sigma=\sigma(\lambda)=\lambda^{\frac{1}{2}} e^{\frac{i \theta}{2}}+\lambda^{-\frac{1}{2}} e^{-\frac{i \theta}{2}}.\]

Since \(\sigma^2=\tau+2\) we obtain that \(\sigma\) has degree at most 2 over \(\Q(\tau)\). It will be useful for us that \(\Q(\sigma)=\Q(\tau)\). For this, we only need to prove that \(\sigma\) has degree 2 over \(\Q\) which we will turn to prove below. 

\begin{lemma}\label{lem:factoring p(x^2)}
Let \(\ell >0\). If $\lambda=e^\ell \in \Sal_2$ has complex conjugates \(e^{\pm i \theta}\), then  we can choose \(g(x)\) in Definition \ref{def:square-rootable} in the form $$g(x)=(x-e^{\frac{\ell}{2}})(x-e^{-\frac{\ell}{2}})(x\pm e^{\frac{i\theta}{2}})(x\pm e^{-\frac{i\theta}{2}}).$$  
\end{lemma}

\begin{proof}
Let $f(x)$ and $g(x)$ as in \Cref{def:square-rootable}, with $K=\Q$. The roots of $f$ are $e^{\pm \ell}, e^{\pm i\theta}$, then the roots of $f(x^2)$ over $\C$ are all simple roots and they are precisely $\pm e^{\pm\frac{\ell}{2}}, \pm e^{\pm\frac{i\theta}{2}}$, and $f(x^2)=g(x)g(-x)$. By interchanging $g(x)$ and $g(-x)$ if necessary, we can assume that $e^{\frac{\ell}{2}}$ is a real root of $g(x)$. Since $g(x)$ is palindromic, $e^{-\frac{\ell}{2}}$ is also a root of $g(x)$. This implies that $-e^{\frac{\ell}{2}}$ and $-e^{-\frac{\ell}{2}}$ are roots of $g(-x)$ and then the other two roots of $g(x)$ come in a complex conjugate pair, that is, either the pair $e^{\frac{i\theta}{2}},e^{-\frac{i\theta}{2}}$ or the pair $-e^{\frac{i\theta}{2}},-e^{-\frac{i\theta}{2}}$.
\end{proof}

\begin{proposition}\label{prop:determinant_square_rootable}
If \(\lambda \in \Sal_2\), then \(\sigma(\lambda)\) is an imaginary quadratic integer. Moreover, \(\Q(\sigma(\lambda))=\Q(\sqrt{d_\lambda})\).  
\end{proposition}
\begin{proof}
First we note that \(\sigma(\lambda)^2=\tau(\lambda)+2\), which implies that \(\Q(\tau)\subset \Q(\sigma)\). By \Cref{pro:field of tau} we only need to check that \(\sigma(\lambda)\) is quadratic over \(\Q\). By \Cref{lem:factoring p(x^2)}, we can factorize \(F_\lambda(x^2)=g(x)g(-x)\) where the roots of \(g\) are
\(\{\lambda^{\frac{1}{2}},\lambda^{-\frac{1}{2}},\zeta,\zeta\inv\}\), with \(\zeta = \pm e^{\frac{i \theta }{2}}\). By definition, we write \(g(x)=X^4+\rho X^3+lX^2+\rho X+1\) with \(\rho \in \R\) and \(l \in \Z.\)

We note that \(\sigma(\lambda)= \pm (\lambda^{\frac{1}{2}} \zeta+\lambda^{-\frac{1}{2}} \zeta\inv)\), and hence we can assume that \(\zeta=e^{\frac{i\theta}{2}}\). Arguing as in the proof of the  \Cref{pro:field of tau}, we get that \(\sigma(\lambda)\) satisfies the equation \(X^2-(\sigma+\bar{\sigma})X+\sigma\bar{\sigma}=0,\) with
\begin{align*}
   \sigma+\bar{\sigma}= &~ 2-l \in \Z \\
   \sigma\bar{\sigma}= &~ \lambda+\lambda\inv+\e+\f=-a  \in \Z 
\end{align*}

Hence, \(\sigma(\lambda)\) is a quadratic integer.

\end{proof}

For a Salem number $\lambda=e^{\ell}, \ell>0,$  of degree \(4,\) with conjugates $\lambda^{-1},e^{\pm i\theta}$, we denote by $L_{\lambda}$ the number field $\Q(e^{\frac{\ell+i\theta}{2}})$.

\begin{proposition}\label{prop: Llambda2 vs Llambda}
    Let $\lambda \in \Sal_2$. Then $L_\lambda$ is a number field of degree $4$ containing $\Q(\sqrt{d_\lambda})$, and $L_{\lambda}=L_{\lambda^N}$ for any $N\geq 1$.
\end{proposition}
\begin{proof}
   We first note that $e^{\frac{\ell+i\theta}{2}}$ is a root of the equation \(X^2-\sigma(\lambda)X+1\). By \Cref{prop:determinant_square_rootable}, this implies that $e^{\frac{\ell+i\theta}{2}}$ is a quadratic unit over $\Q(\sqrt{d_\lambda})$. Now, by Dirichlet's Units Theorem, the group of units $\OO_{L_{\lambda}}^*$ of $\OO_{L_{\lambda}}$ has rank 1. This implies that the power of any element in $\OO_{L_{\lambda}}^*$ of infinite order generates $L_\lambda$. In particular, powers of $e^{\frac{\ell+i\theta}{2}}$ generates $L_\lambda$. So, we get that $L_{\lambda}=L_{\lambda^N}$ for any $N\geq 1$.
\end{proof}

The following proposition characterizes the quadratic extensions of imaginary quadratic number fields that have the form $L_\lambda$, for some $\lambda \in \Sal_2$.

\begin{proposition}\label{prop:quadratic extensions of kd}
    Let $k_d=\Q(\sqrt{-d})$ be an imaginary quadratic number field, and $L$ be a quadratic extension of $k_d$ such that \(L \cap \R = \Q\). Then there exists \(\lambda \in \Sal_2\) such that $L=L_\lambda$.

\end{proposition}

\begin{proof}
Let \(L \mid k_d\) be a quadratic extension. Since \(L\) has no real embeddings, the group of units \(\mathcal{O}_L^*\) has rank 1. Let \(\eta\) be a generator of the kernel of the relative norm map \(N_{L \mid k_d}\). Then \(\eta\) has infinite order, \(L=\Q(\eta)\) and \(\eta\inv\) is the image of \(\eta\) under the unique \(k_d\) isomorphism of \(L\).  Since \(L \cap \R=\Q\) by assumption, we can write \(\eta=r\zeta\) with \(r>0\) and \(\zeta \in \mathbb{S}^1\), where \(\zeta\) has infinite order. In particular, \(\zeta\) is not a root of 1.

We will turn to prove that \(r^2\in\Sal_2\) with conjugates \(r^{-2},\zeta^2,\zeta^{-2}\). 

For this, we will reverse the order in the proof of \Cref{prop:determinant_square_rootable}. The non-trivial conjugates of \(\eta\) are \(\eta^{-1},\overline{\eta},\overline{\eta}^{-1}\). If we denote \(\xi = \eta+\eta^{-1} \in \Od\), then 
\begin{align*}
\xi^2 +\overline{\xi}^2 & = \eta^2+\eta^{-2}+\overline{\eta}^2+\overline{\eta}^{-2}+4 = r^2\zeta^2+r^{-2}\zeta^{-2}+r^2\zeta^{-2}+r^{-2}\zeta^{2}+4 \\
\xi  \overline{\xi} & = |\eta|^2+|\eta^{-1}|^2+\eta\overline{\eta}^{-1} +\eta^{-1}\overline{\eta}       =r^2+r^{-2}+\zeta^2+\zeta^{-2}
\end{align*} 

Hence, \(r^2\zeta^2+r^{-2}\zeta^{-2}+r^2\zeta^{-2}+r^{-2}\zeta^{-2}+4\) and \(r^2+r^{-2}+\zeta^2+\zeta^{-2}\) are symmetric functions on the conjugates of \(\eta.\)  Therefore, the polynomial
\[f(x)=(x-r^2)(x-r^{-2})(x-\zeta^2)(x-\zeta^{-2})\]
has rational integer coefficients. Furthermore, let \(g(x)\) be the minimal polynomial of \(\zeta^2\) over $\Q$. Since \(\zeta^2\) is not a root of unit, then \(g\) has degree at least \(4\), which implies that \(g=f\) and hence \(f\) is the minimal polynomial of the Salem number \(r^2\).

Now, we note that if \(t\) and \(s\) are respectively the coefficients of the monomials \(X^3\) and \(X^2\) of the polynomial
\[u(x)=(x-r)(x-r^{-1})(x-\zeta)(x-\zeta^{-1})\]
then 
\begin{align*}
t^2 =& r^2+r^{-2}+\zeta^2+\zeta^{-2}+2(\eta+\eta^{-1}+\overline{\eta}+\overline{\eta}^{-1})+4=\xi\overline{\xi}+2(\xi+\overline{\xi})+4=|\xi+2|^2 \\
s =& 2 + \eta+\eta^{-1}+\overline{\eta}+\overline{\eta}^{-1} = 2+\xi+\overline{\xi},    
 \end{align*}
which implies that \(t^2,s \in \Z\) and $t^2>0$. Writing $u(x)$ as \[u(x)=X^2-\sqrt{t^2}X^3+sX^2-\sqrt{t^2}X+1,\] and noting that \(f(x^2)=u(x)u(-x)\), we conclude that \(r^2\) is square-rootable over \(\Q.\)

\end{proof}

We will finish this section with a proposition which justifies our choice for restricting the field of definition to be \(\Q\). Before that, we need a lemma that relates totally real number fields with square-rootable Salem numbers.

\begin{lemma}\label{lem:tot real vs sqaure rootable}
    A number field $L$ of degree $m$ is totally real if and only if there is $\lambda\in\Sal_m$ such that $L=\Q(\lambda+\lambda^{-1})$. 
\end{lemma}
\begin{proof}
Any totally real number field $L$ of degree $m$ has the form $L=\Q(\beta+\beta\inv)$ for some Salem number $\beta$ of degree $2m$ \cite[Lem. 3.3]{ERT19}. Note that $\Q(\beta+\beta^{-1})=\Q(\beta^2+\beta^{-2}).$ Indeed, $\Q(\beta^2+\beta^{-2})\subseteq\Q(\beta+\beta^{-1})$, and both fields are extensions of $\Q$ of degree \(m\) since $\beta$ and $\beta^2$ are Salem numbers of degree $2m$ \cite[Lem. 7.1]{ERT19}. The result follows taking $\lambda=\beta^2$, which lies in $\Sal_m$ by \cite[Lem. 7.4]{ERT19}.
\end{proof}

\begin{proposition}\label{prop:rational_base_field}

Let \(\mathfrak{C}\) be a commensurability class of arithmetic groups of the first type in \(\Isom(\HH^n)\), and let $m\in\Z$ with $2\leq2m\leq n+1$. If \(\mathfrak{C}\) realizes all Salem numbers in \(\Sal_m\), then the field of definition of \(\mathfrak{C}\) is $\Q$.
\end{proposition}

\begin{proof}
Let \(K\) be the field of definition of \(\mathfrak{C}\) and let $L$ be a totally real number field of degree $m$. By \Cref{lem:tot real vs sqaure rootable}, $L=\Q(\lambda+\lambda\inv)$ for some $\lambda\in\Sal_m$. Since $\mathfrak{C}$ realizes $\lambda$, \Cref{th:ERT 5.2 e 6.3} implies that \(K \subset L\). The result will then follows from the fact that, if \(K\neq \Q\), there exists a totally real number field of degree $m$ over $\Q$ that does not contain $K$. The proof of this uses basic algebraic number theory and we present it for completeness. Suppose that \(K\neq \Q\) and fix \(p\) a prime number with \(p \mid\Delta_K\), where \(\Delta_K\) denotes the discriminant of \(K\). Now, let \(q\) be a prime number such that \(p \neq q\) and \(q \equiv 1 \pmod{2m}\), and consider the cyclotomic extension \(\Q(\zeta_q)\), where \(\zeta_q\) is a \(q-\)th root of 1. We know that \(\Q(\zeta_q+\zeta_q\inv)\) is a totally real number field whose Galois group is isomorphic to the cyclic group of order \(\frac{q-1}{2}\). If \(L\) denotes the unique subfield of \(\Q(\zeta_q+\zeta_q\inv)\) with Galois group of order \(m\), then \(L\) is a totally real number field of degree \(m\) whose discriminant is a power of \(q\). Hence, \(K\) is not a subfield of \(L\).
\end{proof}

\section{Proof of Theorem A}\label{sec:dim 6}

Let us first recall that $\mathcal{C}_4(\Q)$ denotes the Clifford algebra over $K=\Q$ with generators $i_1,i_2,i_3,i_4$ and $i_j^2=-1$ for $j=1,2,3,4$ (see \Cref{sec:vahlen group}). Consider the element 
\[\omega=\frac{1}{2}(-1+i_1+i_2+i_3) \in \mathcal{C}_4(\Q).\]

Put \(x_1=i_1,x_2=i_2,x_3=\omega\) and \(x_4=i_4\) and let $\mathcal{P}_4$ be the power set of $\{1,2,3,4\}$.  For $I=\{\nu_1,\dots,\nu_\kappa\}\in \mathcal{P}_{4}$ with $\nu_1<\dots<\nu_\kappa,$  we set \(x_I=x_{\nu_1}\cdots x_{\nu_k}\), where we adopt the convention $x_{\emptyset}=1$. Let \(\mathcal{Q}\) be the \(\Z\)-module of $\mathcal{C}_4(\Q)$ generated by the elements \(x_I\).

Since \(i_3=1-i_1-i_2+2\omega \in \mathcal{Q},\) we have that 
\(2\mathcal{Q} \subset \mathcal{C}_4(\Z) \subset \mathcal{Q}\), and so \(\mathcal{Q}\) is a lattice of full rank in \(\mathcal{C}_4(\Q).\) Moreover, it follows from the identity \(\omega^2-\omega+1=0\) that \(\mathcal{Q}\) is a ring containing \(1\), i.e. \(\mathcal{Q}\) is an order of \(\mathcal{C}_4(\Q)\). We also note that \(\mathcal{Q}\) is invariant by the anti-involution \(*\)  of \(\mathcal{C}_4(\Q)\) (see \Cref{sec:Clifford algebras}), which enables us to define the group \(\SL_2(\mathcal{Q})\) in the natural way (see \Cref{sec:Clifford algebras}):
 \begin{align*}
        \SL_2(\mathcal{Q})=\left\{ \begin{pmatrix}
            a & b \\ c & d
        \end{pmatrix} : a,b,c,d \in \Gamma_m \cap \mathcal{Q} \cup\{0\},~ ad^*-bc^*=1 \text{ and } ab^*, cd^* \in  \R^{m+1}\right\}.
    \end{align*}

Now, the relation \(2\mathcal{Q} \subset \mathcal{C}_4(\Z) \subset \mathcal{Q}\) implies that the congruence subgroup \(\SL_2(\mathcal{Q})[2]=\mathrm{Ker}(\SL_2(\mathcal{Q}) \to \SL_2(\mathcal{Q}/2\mathcal{Q}))\) is contained in \(\SL_2(\mathcal{C}_4(\Z))\). Since \(\SL_2(\mathcal{C}_4(\Z))\) is contained in \(\SL_2(\mathcal{Q})\), we get that both groups are commensurable. Therefore, by \Cref{th:arithmetic subgroup Vahlen}, \(\SL_2(\mathcal{Q})\) is an arithmetic subgroup of \(\SL_2(\mathcal{C}_4)\). Now, let \(q_6\) be as given in \eqref{eq: canonical form} with $K=\R$, and consider the composition of the morphisms \(\phi\) and \(\psi\)  introduced in \eqref{phi} and \eqref{map psi}. The group 
\begin{equation}\label{eq:gamma}
\Gamma=\phi \circ \psi (\SL_2(\mathcal{Q}))< \SO_{q_6}(\Q)
\end{equation}

thus defines a classical non-cocompact arithmetic subgroup of $\Isom({\HH^6})$. Indeed, it is a well-known consequence of the Hasse-Minkowski Theorem that any indefinite rational quadratic form of rank at least 5 is anisotropic. Hence, by Godement's compactness criterion, any lattice commensurable to \(\SO_{q_6}(\Z)\) is non-cocompact.

We will see that the hyperbolic orbifold \(\mathcal{O}=\Gamma \backslash\HH^6\) satisfies the statement of \Cref{th:dim 6}. By \cite[Proposition 5.3]{EGM87} there exists an isometry between the upper-half and the hyperboloid models of $\HH^6$ which is equivariant with respect to $\phi \circ \psi$. Hence, it is enough to prove that $\SL_2(\mathcal{Q})$ realizes \emph{all} Salem numbers of degree 2 and \emph{all} elements of \(\Sal_2\).

\begin{lemma} \label{lemma: translation length}
Let \(v \in \mathcal{C}_4(\R) \) be a vector satisfying \begin{equation}
    v^2 - tv+s=0
\end{equation}
with \(t^2-4s<0\)  (see \eqref{eq: quadratic equation for vector}), and let \[ g(X)=X^4 - tX^3+(2+s)X^2-tX+1.\] Denote by $\xi$ a largest complex eigenvalue of $g(X)$. Then, the matrix
\begin{equation}
   A=\begin{pmatrix}
0 & -1\\
1 & v
\end{pmatrix}
\end{equation}
is a loxodromic element in \(\SL_2(\mathcal{C}_4)\) with $\ell(A)=2\log|\xi|.$  
\end{lemma}

\begin{proof}
It follows from the definition of the Vahlen group that \(A \in \SL_2(\mathcal{C}_4)\) (see \Cref{sec:vahlen group}). The action of \(A\) in the hyperbolic space \[\HH^6=\{x=x_0+x_1i_1+\ldots+x_5i_5 \in \mathcal{C}_5: x_5>0\}\] is given by
\[A \cdot x= -(x+v)\inv.\]
Thus, the fixed points of \(A\) are the solutions of 
\begin{equation*}
    x=-(x+v)\inv \Longleftrightarrow x^2+vx+1=0.
\end{equation*}

Let \(z, \overline{z}\) be the complex roots of \(X^2-tX+s=0\). Since \(z \notin \R\), we get \(z^2 - 4 \neq 0\) and thus the equation \(u^2+zu + 1=0\) has two distinct solutions \(\xi,\xi\inv\in \C\) where \(|\xi|>1\). We note that \(\R(z)\) and \(\R(v)\) are  \(\R\)-isomorphic to the field \(\R[X]/(X^2-tX+s) \simeq \C\), which implies that there exists a \(\R\)-isomorphism \(\phi:\R(z) \to \R(v)\) with \(\phi(z)=v\). 
Therefore, \(\phi(\xi), \phi(\xi)\inv \in \R(v)\) are the fixed points of \(A\), which proves that $A$ is loxodromic. 

Now, by \cite[Lem. 13 - 15]{Wat93}, there exists \(P \in \SL_2(\mathcal{C}_5)\) such that
\begin{align*}
   P\inv AP =\begin{pmatrix} r(\cos\theta+\sin\theta i_1) & 0 \\ 0 & r^{-1}(\cos\theta - \sin\theta i_1) \end{pmatrix}
\end{align*}
for some $\theta\in [0,\pi]$, and \(2\log(r)\) is equal to the translation length of \(A\). Since $i_1^2=-1$, the matrix \(P\inv AP\) can be regarded as an element in \(M_2(\C)\). Consider the matrix \[B=\begin{pmatrix}  0 & -1\\ 1 & z \end{pmatrix} \in M_2(\C).\] It follows from \(B^2-zB+I=0\) that \(B\) is annihilated by the real polynomial \begin{align*}
g(X) & =(X^2-zX+1)(X^2-\overline{z}X+1) \\
& = X^4-(z+\overline{z})X^3+ (2+|z|^2)X^2-(z+\overline{z})X+1 \\
     & = X^4 - tX^3+(2+s)X^2-tX+1.
\end{align*}

Then the complex roots of $g(X)$ are precisely \(\xi,\xi^{-1},\overline{\xi},\overline{\xi^{-1}}\).  We will see that $r=|\xi|$.

The isomorphism \(\phi:\R(z) \to \R(v)\) sends $z$ to $v$, \(\phi|_{\R}=Id\), and \(\phi\) extends naturally (coefficient-wise) to a monomorphism of \(\R\)-algebras \[\Phi:M_2(\C) \to M_2(\mathcal{C}_5)\]

where $\Phi(B)=A$. Since $\R$ is contained in the center of \(M_2(\mathcal{C}_5)\), we obtain that \(g(A)=g(\Phi(B))=\Phi(g(B))=0.\) Moreover, 
\begin{align*}
    0 = P\inv g(A)P &  = g(P\inv AP) \\
      & = \begin{pmatrix} 
        g(r(\cos\phi+\sin\phi i_1)) & 0 \\ 0 & g(r\inv(\cos\phi-\sin\phi i_1)) 
    \end{pmatrix}.
\end{align*}

Therefore, regarded as a complex number, $r(\cos\phi+\sin\phi i_1)$ is a largest root of $g$. Since $r=|r(\cos\phi+\sin\phi i_1)|$ we get that $r=|\xi|$ as we wanted.
\end{proof}

\begin{proof}[Proof of \Cref{th:dim 6}]
We will prove that  \(\mathcal{O}=\Gamma \backslash \HH^6\) realizes all the set of square-rootable Salem numbers of degree \( \le 4\). If \(\lambda\) is a Salem number of degree \(2\), with \(\lambda+\lambda\inv=N \in \Z_{>2}\), then \(\lambda\) is also a Salem number of degree \(2\) and square-rootable over \(\Q\) (see \Cref{eq:2-Salem are SR}). By Lagrange's Four Square Theorem, we can write \[N+2=a^2+b^2+c^2+d^2,\] with \(a,b,c,d \in \Z.\) Let \(v=ai_1+bi_2+ci_3+di_4 \in \mathcal{Q}\). Since \(\tr(v)=0\) and \(\mathrm{Norm}(v)=N+2,\) we have that \(v\) satisfies the equation \(v^2+(N+2)=0\). By \Cref{lemma: translation length}, the matrix \[B=\begin{pmatrix}
    0 & -1 \\ 1 & v
\end{pmatrix}\]
is loxodromic with translation length \(2\log(r),\) where \(r\) is the largest root of 
\[g(X)=X^4-NX^2+1=(X-\lambda^{1/2})(X-\lambda^{-1/2})(X+\lambda^{1/2})(X+\lambda^{-1/2}).\]
Thus, \(\ell(B)=2\log(\lambda^{1/2})=\log(\lambda).\)

Now, let \(\lambda\in\Sal_2\) with conjugates \(\lambda\inv,\e,\f\) and minimal polynomial \(F_\lambda(X)=X^4+aX^3+bX^2+aX+1\). The complex number \(\sigma=\sigma(\lambda)=\lambda^{\frac{1}{2}}e^{\frac{i\theta}{2}}+\lambda^{-\frac{1}{2}}e^{-\frac{i\theta}{2}}\) is an imaginary quadratic integer that satisfies the equation \(\sigma^2-(2-l)\sigma-a=0\) for some \(l \in \Z\) (see \Cref{prop:determinant_square_rootable}).
 We claim that there is a vector $\alpha\in\mathcal{Q}$ such that
\begin{equation}\label{quadratic equation}
    \alpha^2-(2-l)\alpha-a =0
\end{equation}
Let \(\Delta=(2-l)^2+4a\) be the discriminant of  $X^2 - (2-l)X -a$. The complex roots of \eqref{quadratic equation} are precisely $\sigma,\overline{\sigma}$, which are not real, and then \(\Delta<0.\) By Lagrange's Four Square Theorem again, we can write \[-\Delta=p_1^2+p_2^2+p_3^2+p_4^2\] with \(p_i \in \Z,\) for \(i=1,2,3,4.\) Let \[\alpha_0=p_1i_1+p_2i_2+p_3i_3+p_4i_4 \in \mathcal{C}_4(\Z).\]

Then, $\alpha_0^2=-\Delta$, and we consider \[\alpha = \frac{l-2+\alpha_0}{2}.\]

Since \(\alpha\) is a vector with \(\tr(\alpha)=l-2\) and \(\rm{Norm}(\alpha)\)\(=-a\), we get that \(\alpha\) satisfies \Cref{quadratic equation}. To see that $\alpha$ lies in \(\mathcal{Q}\), we divide the analysis in two cases, according to the parity of $l$.
\begin{itemize}
    \item \(l\) is even: In this case, it holds that $-\Delta\equiv 0\mod 4$, and by writing $-\frac{\Delta}{4}$ as sum of four squares integers we can suppose that all the $p_i$ are even, i.e. \(p_i=2p'_i\), with \(p'_i \in \Z\) for each \(i=1,2,3,4\), and then 
\begin{align*}
    \alpha & = \frac{l-2+\alpha_0}{2}\\
    & =\left(\frac{l}{2}-1-p'_3\right)+(p'_1-p'_3)i_1+(p'_2-p'_3)i_2+2p'_3\omega+p'_4 i_4\in\mathcal{Q}.
\end{align*} 

\item \(l\) is odd: In this case we have \(-\Delta \equiv -(l-2)^2 \equiv 3 \mod 4\). Hence, without loss of generality, we can suppose that \(p_i\) is odd for \(i=1,2,3\) and \(p_4\) is even.
Then, we can write 
\begin{align*}
  \frac{l-2}{2}=\frac{1}{2}+l' & \text{ and } ~ \frac{p_i}{2}=\frac{1}{2}+p'_i
\end{align*}  for \(i=1,2,3\) and \(p_4=2p'_4\) with \(l',p'_1,\ldots,p'_4 \in \Z\). In this case we have
\begin{align*}
 \alpha & = \frac{l-2+\alpha_0}{2}     \\
   & = (l'+1-p'_3) + (p'_1-p'_3)i_1 +(p'_2-p'_3)i_2+(1+2p'_3)\omega+p'_4 i_4
\end{align*}

\end{itemize}

Now, by \Cref{lemma: translation length} the matrix
\begin{equation}
    A=\begin{pmatrix}
0 & 1\\
-1 & \alpha
\end{pmatrix}
\end{equation}
is a loxodromic element in \(\SL_2(\mathcal{Q})\), whose translation length is equal to \(\log(\lambda),\) since \(\lambda^{\frac{1}{2}}\) is the largest norm of a root of the polyomial 
\begin{align*}
    g(X)&=X^4-(2-l)X^3+(2-a)X^2-(2-l)X+1\\
        &= (X^2-\sigma X+1)(X^2-\overline{\sigma}X+1) \\
    &=
(X-\lambda^{\frac{1}{2}}e^{\frac{i\theta}{2}})(X-\lambda^{-\frac{1}{2}}e^{-\frac{i\theta}{2}})(X-\lambda^{\frac{1}{2}}e^{-\frac{i\theta}{2}})(X-\lambda^{-\frac{1}{2}}e^{\frac{i\theta}{2}}).
\end{align*}

Therefore, the element \(\gamma =(\phi\circ\psi) (A) \in \Gamma\) is loxodromic with translation length \(\log(\lambda)\), which proves that $\Gamma$ realizes $\lambda$. 

To finish the proof, let now $\lambda$ be any Salem number of degree $4$. By \cite[Lem. 7.4]{ERT19} we get that $\lambda^2\in\Sal_2$, and then there exists a closed geodesic in $\mathcal{O}$ with length equal to $2\log(\lambda)$.

\end{proof}

\begin{corollary}
Let \[q_6(x_0,x_1,x_2,x_3,x_4,x_5,x_6)=x_0^2-x_1^2-x_2^2-x_3^2-x_4^2-x_5^2-x_6^2.\]
There exists a non-cocompact arithmetic lattice \(\Gamma<\Isom^{+}({\HH^6})\) commensurable with \(\SO_{q_6}(\Z)\) that realizes all the set \(\Sal_2\).
\end{corollary}
            
\section{Proof of Theorem B}\label{sec:dim 4-5}

\begin{theorem}
There is no a commensurability class of arithmetic groups of the first type in \(\Isom(\HH^3)\) which realizes all elements of \(\Sal_2\).
\end{theorem}

\begin{proof}
  Suppose, by contradiction, that such class \(\mathfrak{C}\) does exist. By \Cref{prop:rational_base_field} we can suppose that the field of definition of \(\mathfrak{C}\) is \(\Q\). Let $\lambda \in \Sal_2$ with minimal polynomial \(F_\lambda\). If $\Gamma$ realizes $\lambda,$ there is $T\in\Gamma$ with characteristic polynomial $p_T$ such that \(p_T(\lambda)=0\), i.e. \(F_\lambda \mid p_T\) and both polynomials are monic and have the same degree. This implies that $p_T(x)=F_\lambda(x)$, and by \cite[Cor. 9.2]{BF15} $F_\lambda(1)F_\lambda(-1) \equiv \det(q)  \mod \SQ$. Furthermore, it follows from the  \Cref{thm:discriminant vs F(1)F(-1)} that \(d_\lambda \equiv F_\lambda(1)F_\lambda(-1) \mod \SQ.\) Hence, the commensurability class \(\mathfrak{C}\) realizes a Salem number \(\lambda\) only if \(d_\lambda \equiv \det(q) \mod \SQ.\)
    
    If we consider a positive rational prime \(p\) with \(-p \not\equiv \det(q) \mod \SQ \), we can invoke \Cref{prop: Salem of given d} in order to obtain a Salem number \(\lambda_0 \in \Sal_2\) with \(d_{\lambda_0} \equiv -p \mod \SQ,\) which implies that \(\lambda_0\) cannot be realized by any group commensurable to $\Gamma$.  
\end{proof}

We recall that any Salem number of degree $4$ is denoted by $\lambda=e^{\ell}, \ell>0$ with conjugates $\lambda^{-1},e^{\pm i\theta}$. We also denote by $L_{\lambda}$ the number field $\Q(e^{\frac{\ell+i\theta}{2}})$.

\begin{proposition}\label{thm:Salem vs splitting}
Let $q$ be a quadratic form over $\Q$ with signature $(1,4)$ over $\R$, and let $\lambda=e^\ell \in \Sal_2$. If $\SO_q(\Q)$ realizes $\lambda$, then $\Ll$ splits $c(q).$
\end{proposition}

\begin{proof}
Let $T\in \SO_q(\Q)$ with eigenvalue $\lambda$. Then, $T^2$ is conjugated in $\SO_q(\R)\simeq \SO(1,4)$ to the matrix
    \begin{center}
\[A=\begin{pmatrix}
\cosh(2\ell) & \sinh(2\ell) & 0 & 0 & 0\\
\sinh(2\ell) & \cosh(2\ell) & 0 & 0 & 0\\
0 & 0 & \cos(2\theta) & \sin(2\theta) & 0\\
0 & 0 & -\sin(2\theta) & \cos(2\theta) & 0\\
0 & 0 & 0 & 0 & 1
\end{pmatrix}.\]
\end{center}

Since $T^2$ lies in the image of the map   $\varphi:\Spin_q(\Q)\rightarrow \SO_q(\Q),$ which sends $s$ to $\varphi_s(x)= sxs^{-1}$
(see \Cref{eq:index spin vs so}), there exists $s\in\Spin_q(\Q)$ with $\varphi(s)=T^2$. In particular, $s$ is conjugate in $\Spin_q(\R)\simeq\Spin(1,4)$ to the element 
  
   $$t=\left(\cosh{\ell}+f_0f_1\sinh{\ell}\right)\left(\cos{\theta}+i_1i_2\sin{\theta}\right).$$

  Observe that $f_0f_1$ and $i_1i_2$ commute, and then generate a commutative subalgebra of $\mathcal{C}^{+}(q)$. With this, a direct computation shows that $t$ is root of the polynomial
        \begin{equation}\label{eq:min pol of e^{l+itheta}}
          r(x)=x^4-(4\cosh{\ell}\cos{\theta})x^3+(2+4\cosh^2{\ell}-4\sin^2{\theta})x^2-(4\cosh{\ell}\cos{\theta})x+1.
        \end{equation}

 By \Cref{lem:min poly beta}, $r(x)$ is the minimal polynomial of $e^{\ell+i\theta}$ over $\Q$. Since $r(x)\in\Z[x]$, the equation $r(x)=0$ is invariant by conjugation, and then $r(s)=0$. Therefore, the $\Q$-subalgebra $\Q[s]$ of $\mathcal{C}^{+}(q)$ is isomorphic to $\Q[e^{\ell+i\theta}]=\Q(e^{\ell+i\theta})=L_{\lambda^2}$. 
 
 We are assuming that \(\lambda\) is square-rootable. Thus, \Cref{prop: Llambda2 vs Llambda} tells us that \(\Ll=L_{\lambda^2}\) has degree $4$ over $\Q$.
 In other words, $\mathcal{C}^{+}(q)$ contains a field of degree $4$ over $\Q$ isomorphic to $\Ll$. Since $\mathcal{C}^{+}(q)$ has dimension 16 over $\Q$ we have that $\Ll$ splits $\mathcal{C}^{+}(q)$ (see \cite[Thm 4.4]{T09} or \cite[Corollary 3.6]{Milne11}), i.e. $\Ll$ splits $c(q)=[\mathcal{C}^{+}(q)]$.
 
\end{proof}

Now, we recall that \(\Od\) denotes the ring of integers of \(k_d=\Q(\sqrt{-d})\), for a square-free integer \(d>0.\)

\begin{lemma}\label{lem: quad ext 1}
    Let $u\in\Od-(\Od^2\cup\Q).$ Then $\Q(\sqrt{u})=k_d(\sqrt{u})$ and $k_d(\sqrt{u})\cap\R=\Q.$
\end{lemma}

\begin{proof}
    Since $u\notin(\Od^2\cup\Q)$ the minimal polynomial of $\sqrt{u}$ over $\Q$ is given by \[p_{\sqrt{u},\Q}(x)=(x^2-u)(x^2-\bar{u})=x^4-(u+\bar{u})x^2+|u|^2,\]

    and then $\sqrt{u}$ has degree $4$ over $\Q$, so $\Q(\sqrt{u})=k_d(\sqrt{u})$. Suppose now that $\beta\in k_d(\sqrt{u})\cap\R$. Then $\Q(\beta)$ is a real subfield of $k_d(\sqrt{u})$. Denote by $p_{\sqrt{u},\Q(\beta)}(x)$ the minimal polynomial of $\sqrt{u}$ over $\Q(\beta)$, which must be a divisor of $p_{\sqrt{u},\Q}(x)$ over $\Q(\beta)$. However, since $u$ is not real, the only divisor of  $p_{\sqrt{u},\Q}(x)$ with real coefficients is $p_{\sqrt{u},\Q}(x)$ itself, so $p_{\sqrt{u},\Q(\beta)}(x)=p_{\sqrt{u},\Q}(x)$ and then $[k_d(\sqrt{u}):\Q(\beta)]=4$. Therefore, $\Q(\beta)=\Q$ and $\beta\in\Q$. 
\end{proof}

\begin{lemma}\label{lem: quad ext 2}
    Let $\mathfrak{p}$ be a prime ideal in $\gd$. There is  \(\lambda \in \Sal_2\) such that $L_\lambda$ is a quadratic extension of $k_d$ that splits $\mathfrak{p}$. 
\end{lemma}

\begin{proof}
    Let $\mathfrak{q}$ be a prime ideal of $\gd$ different to $\mathfrak{p},$ and which does not divide \(2\). By the Chinese Remainder Theorem, there exists $u\in\gd$ satisfying the congruences
    \begin{align*}
        u&\equiv1\mod\mathfrak{p}^5,\\
        u&\not\equiv\zeta^2\hspace{-2mm}\mod\mathfrak{q},\text{\hspace{2mm}  for any } \zeta\in\gd,\\
        u&\not\equiv 0\mod2.
    \end{align*}

If $u\in\Q$, let $\beta\in(\mathfrak{p}^5\mathfrak{q})\cap(2)-\Q$ and set $u'=u+\beta$. So, we can assume that $u\notin\Q$.  The second congruence condition implies that $u\notin\Od^2$. By \cite[Thm. 118-119]{H91} $\mathfrak{p}$ splits in $L=k_d(\sqrt{u})$. Since $u\in\Od-(\Od^2\cup\Q)$, by \Cref{lem: quad ext 1} and \Cref{prop:quadratic extensions of kd}, $L=L_\lambda$ for some Salem number $\lambda \in \Sal_2$.
\end{proof}

For a quaternion algebra over $\Q$, we recall that $\Ram_f(B)$ denotes the set of finite places of $\Q$ that ramifies $B$. 

\begin{proposition}\label{prop:Ll doesnot split Witt}
    Let $B$ be a quaternion algebra over $\Q$ such that $\Ram_f(B)\neq\emptyset.$ Then, there is \(\lambda \in \Sal_2\) such that $L_\lambda$ ramifies $B$.
\end{proposition}

\begin{proof}
    Let $p$ be a rational prime that ramifies $B$. We choose $d>0$ such that $p$ splits on $k_d$, and we consider $\mathfrak{p}$ a prime divisor of $p\Od$. By \Cref{lem: quad ext 2} there is $\lambda$ such that \(k_d \subset L_\lambda\) and $\mathfrak{p}$ splits at $L_\lambda$. 

 Since $p$ splits in $k_d$, the $k_d$-quaternion algebra $B\otimes_{\Q}k_d$ ramifies at $\mathfrak{p}$ (see \cite[Thm 2.4]{R92}). Similarly, $(B\otimes_{\Q}k_d)\otimes_{k_d}L_\lambda$ ramifies at the prime factors of $\mathfrak{p}$ in $L_\lambda$ as a quaternion algebra over $L_\lambda$. Since \(B\otimes_{\Q}L_{\lambda}=(B\otimes_{\Q}k_d)\otimes_{k_d}L_{\lambda}\) we get that $L_\lambda$ ramifies $B$.
\end{proof}

\begin{theorem}\label{thm:nonexistence dim 4}
There is no a commensurability class of arithmetic groups of $\Isom(\HH^4)$ which realizes all Salem number in $\Sal_2$.
\end{theorem}

\begin{proof}

Suppose the contrary and let $\Gamma<\Isom{(\HH^4)}$ be any representative of a commensurability class which realizes a given Salem number \(\lambda \in \Sal_2.\)
By the classification of arithmetic groups, $\Gamma$ is of the first type and classical \cite[Lemma 4.2]{ERT19}. By  \Cref{prop:rational_base_field}, $\Gamma$ is defined over $\Q$.  Then, there is an admissible rational quadratic form $q$ such that $\Gamma\subset\SO_q(\Q)$ and $\Gamma$ is commensurable with $\SO(q,\Z)$. The Witt invariant \(c(q)\) is an invariant of the commensurability class of $\Gamma$. Denote by $D$ the quaternion algebra over $\Q$ representing $c(q)$ in $\mathrm{Br}(\Q)$. By \Cref{thm:Salem vs splitting}, $D$ splits at $L_\lambda$.

On the other hand, since $D$ ramifies at infinity (\Cref{lem:ram at inf}), there is at least one finite prime $p$ that ramifies $D$. By \Cref{prop:Ll doesnot split Witt} there exists $\lambda\in\Sal_2$ such that $L_\lambda$ ramifies $D$. This contradiction implies the result.
 
\end{proof}

We will see now that \Cref{thm:nonexistence dim 4} also holds in dimension 5. We recall that for a rational quadratic space \((V,q),\) of dimension \(r\), its \emph{signed determinant} is the class in $\SQ$ of the number \(\delta = (-1)^{\frac{r(r-1)}{2}}\det(q)\).

\begin{proposition}\label{thm:Salem vs splitting dim 5}
    Let $q$ be a rational quadratic form  with signature $(1,5),$ and consider \(\delta=-\det(q)>0\). Let $D$ be a quaternion algebra representing $c(q)$ in $\Br(\Q)$.  If $\SO_q(\Q)$ realizes \(\lambda\in\Sal_2\), then $\Ll$ splits $D\otimes_\Q\left(\frac{\delta,-1}{\Q}\right)$. 
\end{proposition}

\begin{proof}
 Suppose that $T\in\SO_q(\Q)$ realizes $\lambda$. Let \(p_T \in \Q[X]\) the characteristic polynomial of \(T\). Since \(p_T(\lambda)=0\) we have that \(p_T(x)=F_\lambda(x)g(x)\) for some \(g(x) \in \Q[X]\) of degree 2, where \(F_\lambda(x)\) is the minimal polynomial of \(\lambda\). By \cite[Thm. 5.2]{ERT19}, \(g\) is a product of cyclotomic polynomials of degree 1 or 2. There exists $N\leq 6$ such that the loxodromic element \(T^N\) has eigenvalue \(1\). By applying \Cref{prop: Llambda2 vs Llambda} we also have \(L_\lambda=L_{\lambda^N}\). Thus, we can assume that $T$ has an eigenvalue equal to 1.
 
 We adapt some arguments as in \cite{M11}. Let $x_1$ be a corresponding eigenvector. Note that $q(x_1)$ is negative, and by rescaling, we can assume that $q(x_1)=-1$. Let $V_1=\langle x_1\rangle$ and $V_0$ with $V=V_0\perp V_1$. Then $\det(V)=-\det(V_0)$ and the signed determinants of $V_0$ and $V$ are equal. Denote by $q_0$ the restriction of $q$ to $V_0$. Then, the quadratic forms $-q$ and $-q_0$ have signature  $(5,1)$ and $(4,1)$ respectively, $-q=-q_0\perp<1>$, and $\SO(q_0,\Q)$ realizes $\lambda$.   Arguing as in the proof of \cite[Thm. 6.1]{M11} for $n=5$
 we get that $c(-q)=c(-q_0)$. By \cite[Eq. (3.16)]{Lam05} we have that 
 \begin{align*}
     c(-q_0)&=c(q_0)\\
     c(-q)&=c(q)\left(\frac{-1,\delta}{\Q}\right),
 \end{align*}
 and then $c(q_0)=[D\otimes_\Q\bigl(\frac{-1, \delta}{\Q}\bigr)]$ in $\Br(\Q)$. By \Cref{thm:Salem vs splitting} we get that $L_\lambda$ splits $c(q_0)$ and the result follows.

\end{proof}

With this, we obtain the analog of \Cref{thm:nonexistence dim 4} for dimension 5.

\begin{theorem}
There is no a commensurability class of arithmetic groups of $\Isom(\HH^5)$ of the first type which realizes all Salem numbers in \(\Sal_2\).
\end{theorem}

\begin{proof}

Suppose that such commensurability class $\mathfrak{C}$ exists. By  \Cref{prop:rational_base_field}, any arithmetic subgroup $\Gamma$ of $\Isom(\HH^5)$ in $\mathfrak{C}$ is defined over $\Q$ and then, there is a rational quadratic form $q$ with signature $(1,5)$ over $\R$, such that $\Gamma$ is commensurable with $\SO(q,\Z)$. Let $\delta>0$ be a representative of the signed determinant of \(q\), and denote by $B$ and $D$ the quaternion algebras such that $[D]=c(q)$ and $[B]=\left[D\otimes_\Q\left(\frac{\delta,-1}{\Q}\right)\right]$ in $\mathrm{Br}(\Q)$. Then, \(B\) and \(D\) do not depend on \(\Gamma\), but only on $\mathfrak{C}$.

Now, note that if \(\lambda\) is realized by some \(\Gamma\) in $\mathfrak{C}$, then some power of \(\lambda\) is realized by \(\SO(q,\Z)\), but \(L_\lambda=L_{\lambda^k}\) for any \(k \geq 1.\) Now, by \Cref{thm:Salem vs splitting dim 5}, $B$ splits at $L_\lambda$ for any $\lambda\in\Sal_2$.

On the other hand, since $D$ ramifies at infinity (\Cref{lem:ram at inf}) and $\delta>0$ we get that $B$ also ramifies at infinity. So, $\mathrm{Ram}_f(B)\neq\emptyset$ and \Cref{prop:Ll doesnot split Witt} implies the existence of some $\lambda$ such that $L_\lambda$ ramifies $B$, and this contradiction concludes the proof.

\end{proof}

\appendix

\section{Discriminant of Salem Numbers}\label{sec:discriminant}

We recall that for an algebraic integer \(\alpha\), the \emph{discriminant} \(d_\alpha\) of \(\alpha\) is defined as the discriminant of its minimal polynomial \(f\), i.e. \[d_\alpha=(-1)^{\frac{n(n-1)}{2}}\prod\limits_{\substack{f(\beta)=0,~ f(\delta)=0\\ \beta \neq \delta}} (\beta - \delta),\]
where \(n\) is the degree of \(f.\) 

It follows from \cite{ERT19}  that we can characterize Salem numbers as algebraic integers with minimal polynomial being realized as the minimal polynomial of isometries of a quadratic spaces of signature \((1,n)\), for some integer \(n>1\). 

The problem of determining when an irreducible polynomial $f$ can be the minimal polynomial of an isometry of some quadratic space goes back to the work of Milnor \cite{M69}. With the advance of this problem, it has appeared the importance of the number \(f(1)f(-1)\) (see \cite{GM02}, \cite{BF15}). In this section we will explore the relation of this number with the discriminant of Salem numbers.
 
As we have mentioned in the Introduction, the following proposition is well-known, but we give a proof for the convenience of the reader.

\begin{proposition}\label{thm:discriminant vs F(1)F(-1)}
Let $m\geq 1$ and $\lambda$ be a Salem number of degree \(2m\)  with minimal polynomial \(F_\lambda\), and discriminant $d_\lambda$. Then
    $d_\lambda=(-1)^{m}F_\lambda(1)F_{\lambda}(-1)$ in $\Q^{\times}/\Q^{\times 2}$.
\end{proposition}

\begin{proof}
If \(\lambda\) is a quadratic algebraic integer with minimal polynomial \(F_\lambda=x^2-Nx+1\), then \(d_{\lambda}=N^2-4\) while \(F_\lambda(-1)F_\lambda(1)=4-N^2\), in this case we have the equality \(d_\lambda=(-1)F_\lambda(-1)F_\lambda(1)\).

Now suppose \(m > 1,\) and consider the number fields $E=\mathbb{Q}(\lambda)$ and $K=\mathbb{Q}(\lambda +\lambda^{-1})$. As vector spaces over $\Q$ we know that $\dim_\Q(E)=2m$ and $\dim_\Q(K)=m$. Let $\iota:E\rightarrow E$ be the unique $K$-isomorphism of \(E\), which sends $\lambda$ to $\lambda^{-1}$. Over $E$ we consider the quadratic form $q^\lambda(x)=\tr_{E \mid \mathbb{Q}}(x \cdot \iota(x))$. The result will follow by comparing the determinant of $q^\lambda$ with that of the trace form $q(x)=\tr_{E|\Q}(x\cdot x)$.

A direct computation shows that $K$ and $(\lambda-\lambda^{-1})K$ are subspaces of $E$ which are orthogonal with respect to both $q^\lambda$ and $q$. Therefore
\begin{align*}
  q^\lambda&=q^\lambda_{|K}\oplus q^\lambda_{|(\lambda-\lambda^{-1})K}\\
  q&=q_{|K}\oplus q_{|(\lambda-\lambda^{-1})K}  
\end{align*}

Since $q^\lambda_{|K}(x)=2\tr_{K|\Q}(x^2)$ and $q^\lambda_{|(\lambda-\lambda^{-1})K}(x)=-2\tr_{K|\Q}((\lambda-\lambda^{-1})^2x^2)$ we get that $$\det(q^\lambda)=(-4)^m\N_{K}((\lambda-\lambda^{-1})^2)\Delta_{K}^2,$$ where $\N_K$ denotes the norm function on $K$, and $\Delta_K$ denotes the discriminant of $K$. Similarly, $q_{|K}(x)=2\tr_{K|\Q}(x^2)$ and $q_{|(\lambda-\lambda^{-1})K}(x)=2\tr_{K|\Q}((\lambda-\lambda^{-1})^2x^2)$, so $$\det(q)=4^m\N_{K}((\lambda-\lambda^{-1})^2)\Delta_{K}^2.$$

This implies that $\det(q)=(-1)^{m}\det(q_{\lambda})$. From the definition of the discriminant of $\lambda$ we know that $d_\lambda=\det(q)$ in $\SQ$. On the other hand, note that the map $\phi_{\lambda}(x)=\lambda x$ is an isometry of $(E,q^\lambda)$ with characteristic polynomial $F_\lambda(x)$. By \cite[Cor. 9.2]{BF15} we have that 
\begin{equation}
\det(q^\lambda)=F_\lambda(1)F_\lambda(-1) \mbox{ in } \SQ, 
\end{equation}
and the result follows.
\end{proof}
    
\begin{proposition}\label{prop: Salem of given d}
  Given a square-free integer $d>0$ and an integer \(m>1\), there exists a Salem number $\lambda$ of degree \(2m\) with $d_\lambda \equiv (-1)^{m+1}d \mod \SQ$. Furthermore, we can always choose \(\lambda\) to be square-rootable over \(\Q\).
\end{proposition}

\begin{proof}
 Consider the rational quadratic form $q_d$ given by 
  \begin{align*}
      q_d(x_0,\ldots,x_{2m-1}) = dx_0^2-x_1^2-\cdots-x_{2m-1}^2 .
  \end{align*}
The quadratic form $q_d$ has been chosen, so $\det(q_d)=-d$ and signature \((1,2m-1).\) Hence, the group $\Gamma=\SO(q_d,\Z)$ is a classical arithmetic group of first-type.

Let $T\in\Gamma$ be a loxodromic element with characteristic polynomial $p_T(x) \in \Z[x]$. Assume that $p_T$ is irreducible over $\Q$. Then $p_T$ is the minimal polynomial of a Salem number $\lambda$ of degree $2m$ (see \cite[Thm 5.2]{ERT19}). Since $T$ is an isometry of $q_d$ we have by \Cref{thm:discriminant vs F(1)F(-1)} and \cite[Cor. 9.2]{BF15} that $$-d=\det(q_d) \equiv p_T(1)p_T(-1) \equiv (-1)^{m}d_\lambda\mod \Q^\times/\Q^{\times2}.$$ So, it remains to prove that $\Gamma$ has at least one loxodromic element with irreducible characteristic polynomial. 

It follows from \cite[Theorem 5.2]{ERT19} that $p_T(x)$ is reducible over \(\Q\) if and only if \(T\) has an eigenvalue contained in the set of roots of unit. Suppose that a root of unit \(\xi\) of order \(k\) is an eigenvalue of \(T\). Since $p_T(x)\in\Q[x]$, the minimal polynomial \(\Phi_k\) of \(\xi\) over $\Q$ divides \(p_T\), hence \(\phi(k) \leq 2m,\) where \(\phi(k)\) is the degree of \(\Phi_k\). Thus, there exists an upper bound \(N(m)\) for \(k\).

It is well-known that there exists a compact group \(M \subset \SO(q_d,\R)\), isomorphic as a Lie group to \(\SO(2m-2),\) such that any loxodromic element \(T \in \Gamma\) can be conjugated to a product \(A\Theta\), where \(\Theta\) and \(A\) commute, \(\Theta \in M,\) and \(A\) is diagonalizable with eigenvalues \(1,e^\ell,e^{-\ell}\), where \(\ell\) is the translation length of \(T.\) Moreover, \(\Theta\) represents the conjugacy class of the holonomy of the closed geodesic induced by \(T\) in the orbifold \(\Gamma \backslash \HH^{2m-1}\) (more details of this equivalence can be found in the Introduction of \cite{SW99} or in \cite[Section 2]{MMO14}). In this way, we conclude that if \(T \in \Gamma\) is loxodromic, and its representative of the holonomy does not have eigenvalues of finite order, then the characteristic polynomial of \(T\) is irreducible.

Now, we note that for any \(\zeta \in \mathbb{S}^1,\) the set of transformations in \(\SO(2m-2)\) which have eigenvalue \(\zeta\) is determined by the zero set of the real analytic map \(X \to \det(X^2-(\zeta+\bar{\zeta})X+\mathrm{I})\). Since \(2m-2\) is even, this map is not constant, and then it is a compact set of null Haar measure. Let \(\Sigma\) be the set of operators in \(\SO(2m-2)\) that have an eigenvalue of order bounded by \(N(m)\). We conclude that $\Omega=\Sigma^{c}$ is an open set, invariant by conjugation, and its boundary has zero measure. This allows us to apply the equidistribution of holonomies to $\Omega$ \cite[Corollary of Theorem 1]{SW99} in order to guarantee the existence of many loxodromic elements in \(\Gamma\), whose eigenvalues have infinite order. Hence, the characteristic polynomials of these elements are irreducible as desired. 

If \(T \in \SO(q_d,\Z)\) realizes the Salem number \(\lambda\), then \(T^2\) realizes \(\lambda^2\), and by \cite[Lemma 7.1 and 7.4]{ERT19} \(\lambda^2\) is a Salem number with the same degree that \(\lambda\) and it is square-rootable over \(\Q.\) 

\end{proof}

\section*{Acknowledgments}

We wish to thank Artūras Dubickas for bringing to our attention the reference \cite{D16}, which contains another proof of Proposition A. We are also indebted to Eva Bayer-Fluckiger and Mikhail Belolipetsky for their valuable suggestions and insights.
The author C. D\'oria would like to thank the Isaac Newton Institute for Mathematical Sciences, Cambridge, for support and hospitality during the programme Operators, Graphs, Groups, where work on this paper was undertaken. This work was supported by EPSRC grant EP/Z000580/1, CNPq grant 408834/2023-4 and FAPITEC/SE Grant 019203.01303/2024-1.


\bibliographystyle{plain}
\bibliography{bib2}
\end{document}